\documentclass{article}    
  \usepackage{amsmath,amssymb}
  \newtheorem{definition}{Definition}
  \newtheorem{prop}{Proposition}
  \newtheorem{lem}{Lemma
}
  \newtheorem{remark}{Remark}

    \def\Rset{\mathbb{R}}
    
    \def\Zset{\mathbb{Z}}
    \def\proof{\noindent\mbox{\bf{Proof.}\ }}
    \def\qed{\ \hfill \mbox{$\Box$}}
    \def\rank{\mbox{$\rm{rank}\,$}}
  
\begin{document}

  \title{Controllability to the origin implies 
    state-feedback stabilizability for discrete-time nonlinear systems
    \thanks{
      This work was supported by the Japan Society for the Promotion of 
      Science under Grant-in-Aid for Scientific Research (C) 23560535. 
      This document is the accepted version of the manuscript published in Automatica 
      in Volume 76 with 
      DOI: https://doi.org/10.1016/j.automatica.2016.09.046
      and is a replacement of the preparatory version given in arXiv,
      which the author uploaded at the time instant of its submission to Automatica.  
      If possible, please use the official version instead of this file.  
      $\copyright$2018. This manuscript version is made available under the CC-BY-NC-ND 4.0 license 
      http://creativecommons.org/licenses/by-nc-nd/4.0/       
    }
  }
  \author{Shigeru Hanba
    \thanks{
      Department of Electrical and Electronics Engineering, University of the Ryukyus,
      1 Senbaru Nishihara, Nakagami-gun, Okinawa 903-0213, Japan; email: sh@gargoyle.eee.u-ryukyu.ac.jp
    }
  }
  \maketitle
  
  \begin{abstract}
  The problem of finite-time state-feedback stabilizability of discrete-time nonlinear systems
  has been considered in this technical communique.  Two assertions have been proved.  
  First, if the system is $N$-step controllable to the origin, 
  then there is a state feedback control law for which the trajectory of 
  the closed-loop system converges to the origin in $N$ steps.  
  Second, if the system is asymptotically controllable to the origin and 
  satisfies the controllability rank condition at the origin, 
  then there is a  state feedback control law for which 
  the trajectory of the closed-loop system converges to the origin in finite steps.  
  \end{abstract}
  
  \section{Introduction}
  From control theoretic point of view, one of the most important properties of a system is stability.  
  Controllability assures the existence of an open-loop control law, 
  but in many cases, a state-feedback control law is preferable.  

  For continuous-time systems, the relation between asymptotic controllability and 
  state-feedback stabilizability has been established in \cite{Clarke1997}.  
  In \cite{Clarke1997}, it has been shown that there is a discontinuous state feedback
  stabilizing control law for a system which is asymptotically controllable.  
  Because the discontinuity of the control law arises naturally in stabilization 
  and optimization problems, discontinuous control laws have been studied in many 
  papers, e.g. \cite{Clarke2000,Ceragioli2002,Rifford2002,Clarke2011}.  
  Regarding the time-scale, the stability property is classified into two categories:
  asymptotic property and finite-time property , and recently, there has been 
  a growing interest on the latter property\cite{Haimo1986,Bhat2000,Huang2005,Moulay2006}.  

  For discrete-time systems, problems related to controllability have been extensively 
  studied in \cite{Jakubczyk1990,Albertini1994,Sontag1998}. However, these works do not deal 
  with state-feedback stabilization problem.  
  State feedback stabilization problem of discrete-time nonlinear systems has 
  been studied for past decades (e.g. \cite{Byrnes1993,Kazakos1994,Jiang2001,Jiang2004,Orenelas-Tollez2014}),
  but researches dealing with nonsmooth or discontinuous control laws are relatively rare
  \cite{Meadows1995,Simoes1996}.  
  The connection between controllability and stabilizability has been
  analyzed in an early work by Sontag for piecewise linear systems\cite{Sontag1981}.  
  In \cite{Sontag1981}, both finite-time stability and asymptotic stability 
  has been analyzed, but the scope is limited to piecewise linear systems.  
  More general systems have been deal with in \cite{Kreisselmeier1994,Kellet2004},
  but they concentrate on asymptotic properties.  
  To the best of the author's knowledge,
  the connection between controllability to the origin
  and finite-time state-feedback stabilizability has not been 
  investigated for systems more general than piecewise linear systems.  
  The objective of this technical communique is to fill the gap.  
  
  In the following, we prove two facts.
  First, if a discrete-time nonlinear system is $N$-step controllable to the origin
  (precise definition of this notion is given below), 
  then there is a (possibly discontinuous) state feedback control law for which 
  the trajectory of the closed-loop system converges to the origin in $N$ steps.  
  Second, if the system is asymptotically controllable to the origin and 
  satisfies the controllability rank condition at the origin 
  (again, precise definitions of these notions are given below), 
  then there is a  (possibly discontinuous) state feedback control law for which 
  the trajectory of the closed-loop system converges to the origin in finite steps
  (the required steps may differ for different initial conditions.)  
  Our construction explicitly uses the axiom of choice.  

  A preliminary version of this manuscript is available in arxiv.org
  \cite{Hanba2015}.  

  \section{Main Results}
  Consider a discrete-time time-invariant nonlinear system of the form 
  \begin{equation}
    x(t+1)=f(x(t),u(t)), 
    \label{eq:NLS}
  \end{equation}
  where $x(t) \in \Rset^{n}$ is the state, $u(t) \in \Rset^{m}$ is the control input,
  $t \in \Zset_{\geq 0}$ is the time, and $\Zset_{\geq 0}$ denotes the set of nonnegative integers.  
  It is assumed that $(x,u)=(0,0)$ is an equilibrium of (\ref{eq:NLS}), that is, $f(0,0)=0$.  
  Henceforth, we use the following notations: 
  $u[t_0,t_1]$ denotes the finite sequence of inputs $(u(t_0),\ldots,u(t_1))$, 
  and $u[t_0,\infty)$ denotes the infinite sequence of inputs $(u(t_0),u(t_0+1),\dots)$.
  We identify $u[t_i,t_i]$ with $u(t_i)$, and $u[t_i,t_j]$ with the empty sequence if $t_j<t_i$. 
  For $t \geq t_0$, 
  $\phi(t,t_0,x_0,u)$ denotes the trajectory of (\ref{eq:NLS}) initialized at $t=t_0$
  by $x_0$ and driven by the input (with $\phi(t_0,t_0,x_0,u)=x_0$), 
  and it is also interpreted as the composition of functions defined recursively: 
  $\phi(t,t_0,x_0,u)=f(\phi(t-1,t_0,x_0,u),u(t-1))$. 
  In the subsequent analysis, we sometimes fix $t$ at some $N$ 
  and regard $\phi(N,t_0,x_0,u)$ as a function of $x_0$ and $u[t_0,N-1]$.  
  The notation $\phi(t,t_0,x_0,u)$ implicitly 
  assumes that the input is at least defined over the time interval $[t_0,t-1]$.  
  Because we are dealing with a time-invariant system, 
  the value of the initial time $t_0$ is immaterial.
  Hence, without loss of generality, we assume that $t_0=0$.  

  Properties related to controllability has been used by many researchers for different meanings.  
  The following two definitions are analogous to the one employed in \cite{Sontag1981}.  

  \begin{definition}\label{def:NCO}
    The system (\ref{eq:NLS}) is said to be $N$-step controllability to the origin if
    $\exists N >0$, $\forall x_0 \in \Rset^{n}$, $\exists u[0,N-1]$, $\phi(N,0,x_0;u)=0$.  
  \end{definition}

  \begin{definition}\label{def:ACO}
    The system (\ref{eq:NLS}) is said to be asymptotically controllability to the origin if 
    $\forall x_0 \in \Rset^{n}$, $\exists u[0,\infty)$, $\phi(t,0,x_0;u)$ converges to the origin as $t\rightarrow \infty$.  
  \end{definition}
  
  Another popular notion on controllability is the generalization of
  controllability rank condition of linear systems to nonlinear systems
  given below.  
  \begin{definition}\label{def:RCO}
    The system (\ref{eq:NLS}) is said to be rank controllable at the origin if 
    $\exists N>0$, 
    \begin{equation}
      \rank \frac{\partial \phi(N,0,x_0,u)}{\partial u[0,N-1]}=n
      \label{eq:NL_rank_C}
    \end{equation}
    on an open neighborhood of $(x,u(0),\ldots,u(N-1))=(0,0,\ldots,0)$, \\
    where $\partial \phi(N,0,x_0,u)/ \partial u[0,N-1]$ denotes 
    the partial derivative of $\phi$ with respect to the variable
    $u[0,\ldots,N-1]$ with the re-interpretation that 
    $u[0,\ldots,N-1]$ is the vector $(u^T(0),u^T(1),\ldots,u^T(N-1))^T$,
    and $\cdot^T$ denotes the transpose of a vector.  
  \end{definition}
  In \cite{Sontag1998}, the rank condition (\ref{eq:NL_rank_C}) 
  does not have an explicit name,
  but a tuple $(x,u(0),\ldots,u(N-1))$ that satisfies (\ref{eq:NL_rank_C}) 
  is called `regular.'  Note that we are specializing at 
  $(x,u(0),\ldots,u(N-1))=(0,0,\ldots,0)$.  

  \begin{remark}
    Definition~\ref{def:NCO} and Definition~\ref{def:ACO} are conditions that 
    are not checkable for general nonlinear systems and should blindly be assumed,
    except for polynomial systems with rational coefficients\cite{Nesic1998}.  
    On the other hand, for a discrete-time nonlinear system, 
    $\phi(t,0,x_0;u)$ is merely a composition of known functions,
    hence Definition~\ref{def:RCO} is checkable, although computationally demanding.  
  \end{remark}

  Because we are dealing with a time-invariant system,
  it is preferable that our controller is time-invariant as well.  
  \begin{definition}
    The system (\ref{eq:NLS}) is said to be (globally) finite-time stabilizable by a state feedback controller
    if there is a control law $\upsilon(x)$ defined in $\Rset^{n}$ with the property that
    \[
    \forall x_0 \in \Rset^{n}, \exists N_x>0, \phi(N_x,0,x_0;\upsilon(x))=0.  
    \]
    On the other hand, if
    \[
    \exists N>0, \forall x_0 \in \Rset^{n}, \phi(N,0,x_0;\upsilon(x))=0,
    \] 
    then the system (\ref{eq:NLS}) is said to be (globally) $N$-step stabilizable by a state feedback controller
  \end{definition}
  Because we do not deal with local property in this technical communique,
  henceforth, we omit the term `global.'  

  Our first objective is to show the following proposition,
  which is a straightforward extension of the result given 
  in \cite{Sontag1981} for piecewise linear systems.  
  \begin{lem}\label{lem:NC}
    The system (\ref{eq:NLS}) is $N$-step controllable to the origin
    if and only if it is $N$-step stabilizable by a state feedback controller.  
  \end{lem}
  
  \proof
  One direction is straightforward: 
  if (\ref{eq:NLS}) $N$-step stabilizable by a state feedback controller,
  it is $N$-step controllable to the origin by the input determined by the controller.
  
  For the converse, let $A_0=\{0\}$, and inductively define
  $A_{k}=\{x \in \Rset^n: \exists u, f(x,u) \in A_{k-1}\}$.  
  Because the system is $N$-step controllable to the origin, $A_{N} = \Rset^{n}$.  
  For the sequence $(A_0,A_1,\ldots,A_N)$ and each $x \in \Rset^{n}$, let 
  \begin{equation}
    i_{x}=\min\{i: x \in A_{i} \}.  
    \label{eq:def_i_x}
  \end{equation}
  If $i_{x}=0$, let $\upsilon(x)=0$ (recall that $A_0=\{0\}$).  
  Otherwise, the set $U_{x}=\{u: f(x,u) \in A_{{i_x}-1}\}$ is nonempty.  
  Pick an element $u_x \in U_{x}$ (here, we use the axiom of choice), 
  and let $\upsilon(x)=u_x$.  Because $i_{x}$ is uniquely determined for each $x$,
  $\upsilon(x)$ is well defined.  Let $(x_{t})_{t \in \Zset_{\geq 0}}$ be the trajectory of (\ref{eq:NLS}) initialized with $x_0$
  and driven by the control law $\upsilon(x)$, that is, $x_t=\phi(t,0,x_0;\upsilon(x))$.  
  We prove that $i_{x_t} =0$ for some $t \leq N$ by contradiction.  
  Suppose that $\forall t$, $i_{x_t}>0$.  
  Then, $x_0 \in A_{i_{x_0}}$ for some $i_{x_0} \leq N$ and hence 
  $x_1 =f(x_0,\upsilon(x_0)) \in A_{i_{x_0}-1}$, and by (\ref{eq:def_i_x}),
  $i_{x_1} \leq i_{x_0}-1$.  
  Inductively, assume that $x_j \in A_{i_{x_j}}$ with $i_{x_j} \leq i_{x_0}-j$.  
  Then, $x_{j+1}=f(x_j,\upsilon(x_j)) \in A_{i_{x_j}-1}$, hence 
  $i_{x_{j+1}} \leq i_{x_0}-(j+1)$.  Therefore, for any $j$, $i_{x_j} \leq i_{x_0}-j$.  
  But this is impossible, because $0 \leq i_{x_0} \leq N$ and $0 \leq i_{x_j} \leq N$.  
  \qed

  \begin{remark}
    In \cite{Kellet2004}, the following fact is shown (Theorem 15 of \cite{Kellet2004})
    (in this manuscript, 
    the input space ${\mathcal U}$ 
    and 
    the target set ${\mathcal A}$ of \cite{Kellet2004} are
    are $\Rset^{m}$ and $\{0\}$, 
    respectively, 
    hence we omit them in quoting the result of \cite{Kellet2004}).  
    \begin{quote}
      Let $\sigma: [0,\infty)\rightarrow [0,\infty)$ be a nondecreasing function,
      $\overline{B}_u=\{ u \in \Rset^{m}: \|u\| \leq 1 \}$, and define the set-valued map
      $F(x)$ by $F(x)=f(x,\sigma(\|x\|)\overline{B}_u)$.  If $F(x)$ is continuous in the sense of 
      \cite{Kellet2004}, $F(x)$ is a non-empty compact set for each $x$,
      and (\ref{eq:NLS}) is 
      uniformly globally asymptotically controllable to the origin with $\sigma$ controls
      in the sense of \cite{Kellet2004},
      then there is a feedback function such that the origin is robustly globally asymptotically stable.  
    \end{quote}
    It is to be noted that Lemma~\ref{lem:NC} is not a consequence of Theorem~15 of \cite{Kellet2004}, 
    because we have not required above $F(x)$ to be continuous.  
    As a example, consider the following difference equation defined on the real line: 
    \begin{equation}
      x(t+1)=
      \begin{cases}
        \pi \left \lfloor  \frac{K}{(x(t))^2}  \right \rfloor + \sin (\pi x(t)) u(t), & x(t) \neq 0,\\
        0, & x(t)=0,
      \end{cases}
      \label{eq:CA01}
    \end{equation}
    where $\lfloor r  \rfloor$ denotes the greatest integer not exceeding $r$
    and $K \geq 1$ is a constant.  
    The system (\ref{eq:CA01}) is 2-step reachable to the origin.  
    For, let
    \[
      u(x)=
      \begin{cases}
        -\frac{\pi \left \lfloor \frac{K}{x^2} \right \rfloor}{\sin(\pi x)}, & x \not \in \Zset, \\
        0, & x \in \Zset, 
      \end{cases}
    \]
    where $\Zset$ denotes the set of all integers.  
    If $x_0 \not \in \Zset$, $\phi(1,0,x_0;u)=0$, hence $\phi(2,0,x_0;u)=0$ irrespective of $u(1)$.  
    For $x_0 \in \Zset$ with $x_0^2 > K$, $\phi(1,0,x_0;u)=0$ as well, due to the  effect of $\lfloor \cdot \rfloor$,
    hence $\phi(2,0,x_0;u)=0$.  
    For $x_0 \in \Zset$ with $x_0^2 \leq K$, $\phi(1,0,x_0;u)$ is not an integer, and hence $\phi(2,0,x_0;u)=0$.  

    This result does not follow from Theorem~15 of \cite{Kellet2004}, because the function $F(x)$ obtained from 
    (\ref{eq:CA01}) is not continuous.  For, let $B_u=\{ u \in \Rset^{m}: \|u\| < 1 \}$,
    and $B(x,\delta)=\{z \in \Rset^{n}: \|z-x\| < \delta \{$.  
    The definition of continuity of $F(x)$ in \cite{Kellet2004} requires that
    $\forall x$, $\forall \varepsilon>0$, $\exists \delta>0$, $F(B(x,\delta)) \subset F(x)+\varepsilon B_u$
    (Definition~1 of \cite{Kellet2004}).  In our example, $F(0)=0$, but $F(x)$ is unbounded on the neighborhood of the origin.  
    Hence, our result on stabilizability does not follow from \cite{Kellet2004}.  
  \end{remark}

  \begin{remark}
    Nonlinear finite-horizon model predictive control (MPC) typically assumes that the system is 
    $N$-step controllable to the origin, and shows that the MPC-based state-feedback controller
    asymptotically stabilizes the origin under some conditions on the stage cost
    \cite{Meadows1995,Maciejowski2002}.  Contrary, Lemma~1 assures finite-step stability.  
    It is also to be noted that the controller of Lemma~1 may be of zero robustness, due to its discontinuity.  
  \end{remark}

  Thus far, we have not assumed neither differentiability nor continuity of $f(x,u)$
  and constructed a possibly discontinuous state feedback control law merely under the assumption of 
  $N$-step controllability to the origin.  

  Next, we assume that $f(x,u)$ is at least $C^1$, and show a stronger result.  
  \begin{prop}\label{prop:AC}
    Assume that $f(x,u)$ is $C^1$.
    If the system (\ref{eq:NLS}) is asymptotically controllable to the origin
    and rank controllable at the origin, then it is 
    finite-time stabilizable by a state feedback controller. 
  \end{prop}
  
  \proof
  As in the proof of Lemma~\ref{lem:NC}, 
  let $A_0=\{0\}$, and inductively define
  $A_{k}=\{x \in \Rset^n: \exists u, f(x,u) \in A_{k-1}\}$.  
  We first show that $A_N$ contains an open neighborhood of $x=0$.  
  This is a direct consequence of the implicit function theorem. 
  For, let $p$ be the collection of components of $u[0,N-1]$ for which 
  $\rank \partial \phi(N,0,x,u[0,N-1])/\partial p = n$.  
  Let us rewrite $u[0,N-1]=(p^T,q^T)^T$ by rearranging the variables,
  and rewrite $\phi(N,0,x,u[0,N-1])$ as $\psi(x,p,q)$.  
  Because $\psi(0,0,0)=0$ and $\rank \partial \psi/\partial p =n$,
  by applying the implicit function theorem, 
  $0=\psi(x,p,q)$ may be solved for $p$ in the form $p=h(x,q)$ for some smooth function $h(x,q)$ 
  on a neighborhood of $(x,q)=(0,0)$, and hence
  the trajectory initialized at this neighborhood of $x$ may be driven
  to the origin by applying the corresponding input sequence.  
  Hence $A_N$ contains an open neighborhood of $x=0$ 
  (similar idea to above analysis has been used in Lemma~2 of \cite{Hanba2009}.  )
  Next, we show that $\cup_{k=0}^{\infty}A_k = \Rset^{n}$.  
  For, suppose that $\Rset^{n}\setminus \cup_{k=0}^{\infty}A_k \neq \emptyset$.  
  Let $x_0 \in \Rset^{n}\setminus \cup_{k=0}^{\infty}A_k$.  
  Then, $\phi(t,0,x_0,u) \not \in A_k$ for any $t$, $k$ and $u$ because otherwise $x_0$ is in some $A_k$.  
  However, $A_N$ contains an open neighborhood of the origin.  
  Let the neighborhood be $G$.  Because the system is asymptotically controllable to the origin, 
  by choosing an adequate $u[0,\infty)$, $\phi(t,0,x_0,u) \in G$ residually, whereas
  we have seen that $\phi(t,0,x_0,u) \not \in G$.  This is a contradiction.  
  Hence, $\Rset^{n}\setminus \cup_{k=0}^{\infty}A_k$ must be empty, whence $\cup_{k=0}^{\infty}A_k = \Rset^{n}$.  
  
  Because $\cup_{k=0}^{\infty}A_k = \Rset^{n}$, the control law may be constructed
  in the same fashion as that of Lemma~\ref{lem:NC}.  
  Again, for each $x \in \Rset^{n}$, let 
  $i_{x}=\min\{i: x \in A_{i} \}$.  
  Because $\cup_{k=0}^{\infty}A_k = \Rset^{n}$, $i_{x}$ is well defined and finite.  
  If $i_{x}=0$, let $\upsilon(x)=0$. 
  Otherwise, the set $U_{x}=\{u: f(x,u) \in A_{{i_x}-1}\}$ is nonempty.  
  Pick an element $u_x \in U_{x}$ (here, we use the axiom of choice), 
  and let $\upsilon(x)=u_x$.  
  The proof that this state feedback control law makes the trajectory converge to the origin 
  in finite steps is identical to that of Lemma~\ref{lem:NC}.  
  Note that, for each $x \in \Rset^{n}$, the trajectory initialized at $x$ converges 
  to the origin in at most $i_{x}$ steps.  
  The required steps may differ for different initial conditions.  
  \qed

  \begin{remark}
    In \cite{Hanba2009}, it has been proved that, if the set of permissible initial conditions (denoted by $\Omega$)
    is compact and $\forall x_0 \in \Omega$, 
    $\exists N_{x_0}$, $\exists u[0,N_{x_0}-1]$, 
    i) $\phi(N_{x_0},0,x_0,u)=0$
    and ii) $\rank (\partial \phi(N_{x_0},0,x_0,u)/\partial u[0,N_{x_0}-1])=n$ at $(x_0,u[0,N_{x_0}-1])$, 
    then $\exists N<\infty$, $\forall x_0 \in \Omega$, $N_{x_0} \leq N$, and 
    it is possible to construct a robust stabilizing controller for 
    (\ref{eq:NLS}) through min-max MPC.  
    Proposition~\ref{prop:AC} assures that condition i) of above result may be weakened into 
    asymptotic controllability to the origin.  
  \end{remark}
  
  \section{Conclusion}
  In this technical communique, we have proved two facts.  
  First, a discrete-time nonlinear system is $N$-step controllable to the origin
  if and only if it is $N$-step stabilizable by a state feedback controller.  
  Second, if the system is asymptotically controllable to the origin and 
  satisfies the controllability rank condition at the origin, 
  it is finite-time stabilizable by a state feedback controller.  
  It is also to be noted that the controller may be highly discontinuous 
  and may not be robust.

  \end{document}